\documentclass[12pt]{article}
\usepackage{amsmath}
\usepackage{amssymb}
\usepackage{mathptmx}
\usepackage{amsfonts}
\usepackage[mathscr]{eucal}
\usepackage[dvips]{graphicx}
\usepackage{color}

\newcommand{\mbf}[1]{\ensuremath{\mathbf{#1}}}
\newcommand{\ms}[1]{\ensuremath{\mathscr{#1}}}

\newcommand{\Dir}{\textrm{D}}

\newcommand{\tens}{\otimes}

\newcommand{\restr}{\upharpoonright}

\newtheorem{theorem}{Theorem}

\newtheorem{lemma}[theorem]{Lemma}

\newtheorem{remark}{Remark}

\begin{document}

\begin{center}
{ \Large Essential Self-Adjointness of Anti-Commutative Operators} \\

$\;$ \\
{\large Toshimitsu TAKAESU } \\
\textit{Faculty of Science and Engineering, Gunma University,\\ Gunma, 376-8515, Japan }
\end{center}
\begin{quote}
\textbf{Abstract}.
  In this article,  the self-adjoint extensions of symmetric operators  satisfying  anti-commutation relations are considered. It is proven that an anti-commutative type  of the  Glimm-Jaffe-Nelson commutator theorem follows. Its application to an abstract Dirac operator is  also considered.  
\\
\end{quote}

\section{Introduction and  Main Theorem}
In this article we consider the essential self-adjointness of symmetric  operators satisfying anti-commutation relations. 
Let $H$ be a symmetric operator on a  Hilbert space $\ms{H}$, i.e., $H$ satisfies $H \subset H^{\ast}$. It is said that  $H$ is  self-adjoint if $H = H^{\ast}$ and $H$ is essentially self-adjoint if its closure $\overline{H}$ is self-adjoint.  We are interested in conditions under which a symmetric operator is essentially self-adjoint.   The Glimm-Jaffe-Nelson commutator theorem (e.g., \cite{Ar}; Theorem 2.32,  \cite{RS2}; Theorem X.36) is one criterion for the essential self-adjointness of commutative symmetric operators. 
The  commutator theorem shows that  if a symmetric operator $H$  and a self-adjoint operator $X$ obey a commutation relation on a dense subspace $\ms{D}$, which is a core of $X$, then $H$ is essentially self-adjoint on $\ms{D}$.
Historically, Glimm-Jaffe \cite{GJ72} and Nelson \cite{Ne72} investigate the commutator theorem for quantum field models. Faris-Lavine \cite{FL74}  apply it  to quantum mechanical models  and Fr\"{ohlich}\cite{Fr77} consider a generalization  of the commutator theorem and prove that a  multiple commutator formula follows.

$\;$ \\
The idea of the proof of the commutator theorem is as follows.
Let  $X$ and $Y$ be symmetric operators on a Hilbert space. 
    Then the real part and the imaginary part of the  inner product $\left( X \Psi , \,  Y \Psi   \right)$ for $ \Psi \in \ms{D} (XY) \cap \ms{D} (YX) $ are expressed by 
\begin{align}
&\text{Re} \left( X \Psi , \,  Y \Psi   \right)  \; 
 = \; \frac{  1}{2} \left(  \Psi , \,  \{ X, Y \}  \Psi   \right) \;     , \label{acc} \\
& \text{Im} \left( X \Psi , \, Y \Psi   \right)  \,  \; 
 = \; \frac{ 1 }{2i} \left(  \Psi , \,  [X, Y ]  \Psi   \right)  \;     ,   \label{ccc}
\end{align}
$\;$ \\
respectively, where  $\{ X , \, Y \} = XY + YX $ and $[X,Y ] = XY - YX $.
In the proof of the  commutator theorem,  the imaginary part (\ref{ccc}) is estimated.  In this article,
 we consider symmetric operators satisfying anti-commutation relations. We  estimate   the real part (\ref{acc}) and obtain   an anti-commutative version of the commutator theorem.
Here we  overview the  commutator theorem.

\newpage
Let $H$ and $X$ be linear operators on $\ms{H}$. 
Assume   following conditions.
\begin{quote}
\textbf{(C.1)}   $H$ is symmetric and $X$ is self-adjoint. \\
\textbf{(C.2)} There exists $\delta_{X} > 0 $ such that for all $\Psi  \in \ms{D} (X)$, 
\[
   \delta_{X} (\Psi , \Psi )  \; \leq \;  \left|    ( \Psi , X \Psi  )\right| .  
\]
\textbf{(C.3)}  $X$  has a core $\ms{D}_0$ satisfying $\ms{D}_0 \subset \ms{D} (H)  $, and there exist  constants $a \geq 0 $ and $b \geq 0 $ such that     for all $\Psi \in \ms{D}_0$,
\begin{equation}
\left\|  H \Psi \|  \; \leq \;  a \| X \Psi \right\|  + b \| \Psi \| . \notag 
\end{equation}
\end{quote}
\textbf{{\large Theorem A} (Glimm-Jaffe-Nelson Commutator Theorem}) \\
Let $H$ and $X$ be  operators satisfying \textbf{(C.1)}-\textbf{(C.3)}.
Suppose \textbf{(i)} or \textbf{(ii)} below : \\
\textbf{(i)} There exists a constant $ c_{1} \geq 0 $ such that for all $\Psi \in \ms{D}_{0}$,
\begin{equation}
\left| \frac{}{}  ( H\Psi , X \Psi  ) \, - \, (X \Psi , H \Psi ) \right| \; \leq \; c_{1} \, \left| ( \Psi , X \Psi ) \right| .  \label{4.1}
\end{equation}
\textbf{(ii)} There exists a constant $c_{2} \geq 0 $ such that for all $\Psi \in \ms{D}_{0}$,
\begin{equation}
c_{2} \, \left| ( \Psi , X \Psi ) \right| \; \leq \;  
\left| \frac{}{}  ( H\Psi , X \Psi  ) \, - \, (X \Psi , H \Psi ) \right| \;  .  \label{4.2}
\end{equation}
Then $H$ is essentially self-adjoint on $\ms{D}_{0} $. 

\begin{remark} 
In the commutator theorem, the condition \textbf{(i)} is usually supposed. It is also proven  under the condition \textbf{(ii)} in a similar way to Theorem 1.
\end{remark}

$\; $ \\
The main theorem in this article is as follows.

%%%%%%%%%%%%%%%%%%%%%%%%%%%%%%%%%%%%%%%%%%%%%%%%%%%%%%%%%%%%%%%%%%%%%%%%%%%%%%%%%%%%%%%%%%%%%%%%%%%%%%%%%%%%%%%%%%%%%%%%%%%%%%%%%%%%%%%%%%%%%%%%%%%%%%%%%%%%%%%%%%%%%%%%%%%%%%%%%%%%%%%%%%%%%%%%%%%%%%
\begin{theorem} \label{main}
Assume \textbf{(C.1)}-\textbf{(C.3)}.
In addition suppose that next  \textbf{(I)} or \textbf{(II)} holds. \\
\textbf{(I)} There exists a constant $ d_{1} \geq 0 $ such that for all $\Psi \in \ms{D}_{0}$,
\begin{equation}
\left| \frac{}{}  ( H\Psi , X \Psi  ) \, + \, (X \Psi , H \Psi ) \right| \; \leq \; d_{1} \, \left| ( \Psi , X \Psi ) \right| .  \label{I}
\end{equation}
\textbf{(II)} There exists a constant $d_{2} \geq 0 $ such that for all $\Psi \in \ms{D}_{0}$,
\begin{equation}
d_{2} \, \left| ( \Psi , X \Psi ) \right| \; \leq \;  
\left| \frac{}{}  ( H\Psi , X \Psi  ) \, + \, (X \Psi , H \Psi ) \right| \;  .  \label{II}
\end{equation}
Then $H$ is essentially self-adjoint on $\ms{D}_{0} $. 
\end{theorem}

$\;$ \\
{\large \textbf{(Proof of Theorem 1)}} \\
We  show  that for some  $z \, \in \, \mbf{C} \backslash \mbf{R}$,
$ \text{dim ker} \left( \frac{}{}  (H_{\restr \ms{D}_0 })^\ast
 \, + z^{\sharp} \,  \right) = 0 $ where $z^{\sharp} =z $, $z^{\ast}$.
  Let $\Psi \in \ms{D}  ( (H_{\restr \ms{D}_0 })^\ast) $ and $ \Xi =  X^{-1} \Psi  $.  
Since $ ((H_{\restr \ms{D}_0 })^\ast)^{\ast} = \overline{H_{\restr \ms{D}_0}} $, we have
\begin{equation}
\text{Re} \left( \Xi , \, (  (H_{\restr \ms{D}_0 })^\ast  \, + z^{\sharp} ) \Psi   \right)
\; \; = \; \; \frac{1}{2}\left( \frac{}{} 
 ( \overline{H_{\restr \ms{D}_0}}  \Xi , X \Xi ) \, +  \, ( X \Xi ,  \overline{H_{\restr \ms{D}_0}}  \Xi )  \right)
 \; \; +  \;  \; \text{Re} \, z \, \, (  \Xi , X \Xi )  .  \label{9/25.4}
\end{equation}
First we assume that \textbf{(I)} holds.
Let  $z \in \mbf{C} \backslash \mbf{R} $ satisfying $ |\text{Re} \, z | > \frac{d_{1}}{2} $.  
Since $\ms{D}_0$ is a core of $X$, it follows
 from \textbf{(C.3)} and \textbf{(I)} that 
 $\ms{D} (X) \; \subset  \ms{D} (\overline{H_{\restr \ms{D}_0 }}  )$  and  for all $\Phi \in \ms{D} (X)$, 
\begin{equation}
\left| \frac{}{}  (  \overline{H_{\restr \ms{D}_0 }}\Phi , X \Phi  ) \, + \, (X \Phi ,  \overline{H_{\restr \ms{D}_0 }} \Phi ) \right| \; \leq \; d_{1} \, \left| ( \Phi , X \Phi ) \right| . \label{4.11}
\end{equation}
By (\ref{9/25.4}) and (\ref{4.11}), we have
\begin{equation}
\left| \text{Re} \; \left( \Xi , \, (  (H_{\restr \ms{D}_0 })^\ast  \, + z^{\sharp} ) \Psi   \right) \; \right|
\; \geq \; \left( |\text{Re} \, z | - \frac{d_{1}}{2} \right) | (  \Xi , X \Xi ) | \geq 
  \delta_{X} \left( |\text{Re} \, z | - \frac{d_{1}}{2} \right)  (\Xi , \Xi ) .  \label{9/25.5}
\end{equation}
Since $ \Psi \in $ ker $\left( (H_{\restr \ms{D}_0 })^\ast + z^{\sharp} \right) $, we have
 $  \Xi = X^{-1} \Psi = 0 $ from (\ref{9/25.5}).  Then  we have $\Psi =0$.   Next  we suppose that \textbf{(II)} follows.  Let  $z \in  \mbf{C} \backslash \mbf{R}$ satisfying  $ |\text{Re} \, z | < \frac{d_{2}}{2} $. 
Since $\ms{D}_0$ is a core of $X$, it also  follows
 from \textbf{(C.3)} and \textbf{(II)} that 
 $\ms{D} (X) \; \subset  \ms{D} (\overline{H_{\restr \ms{D}_0 }}  )$  and  for all $\Phi \in \ms{D} (X)$, 
\begin{equation}
 d_{2} \, \left| ( \Phi , X \Phi ) \right|   \; \leq \; 
\left| \frac{}{}  (  \overline{H_{\restr \ms{D}_0 }}\Phi , X \Phi  ) \, + \, (X \Psi ,  \overline{H_{\restr \ms{D}_0 }} \Phi ) \right|  . \label{4.12}
\end{equation}
Then from (\ref{9/25.4}) and (\ref{4.12}), we have
\begin{equation}
\left| \text{Re} \; \left( \Xi , \, (  (H_{\restr \ms{D}_0 })^\ast  \, + z^{\sharp} ) \Psi   \right) \; \right|
\; \geq \; \left(   \frac{d_{2}}{2} - |\text{Re} \, z | \right) | (  \Xi , X \Xi ) | \geq 
  \delta_{X} \left( \frac{d_{2}}{2} - |\text{Re} \, z |  \right)  (\Xi , \Xi ) .  \label{9/25.6}
\end{equation}
Since $ \Psi \in $ ker $\left( (H_{\restr \ms{D}_0 })^\ast + z^{\sharp} \right) $ and $\Xi = X^{-1} \Psi $, we have
   $\Psi =0$    from (\ref{9/25.6}). Thus  the proof is obtained. $\blacksquare$

%%%%%%%%%%%%%%%%%%%%%%%%%%%%%%%%%%%%%%%%%%%%%%%%%%%%%%%%%%%%%%%%%%%%%%%%%%%%%%%%%%%%%%%%%%%%%%%%%%%%%%%%%%%%%%%%%%%%%%%%%%%%%%%%%%%%%%%%%%%%%%%%%%%%%%%%%%%%%%%%%%%%%%%%%%%%%%%%%%%%%%%%%%%%%%%%%%%%%%

\section{Application of Theorem \ref{main}}
We apply  Theorem 1 to a model in supersymmetric quantum mechanics (\cite{Tha}).
Let $\ms{H}$ be a Hilbert space. Let $H$ and $\tau $ be self-adjoint operators on $\ms{H}$.
Assume that $\tau $ is bounded, $\tau^2 = I $ and $\tau \,\ms{D} (H) \subset \ms{D} (H)$. Then $H$ is called an abstract Dirac operator on $\ms{H}$ with unitary involution $\tau$.  We construct an abstract Dirac operator by weakly commuting  operators. Let $X$ and $Y$ be densely defined linear operators on a Hilbert space. The weak commutator of $X$ and $Y$ is defined by  for  $\Phi \in    \cap \ms{D}(X^{\ast}) \cap \ms{D}(Y^{\ast}) $ and for $\Psi \in  \ms{D}(X) \cap \ms{D}(Y)$,
\[
 [X, Y ]^0(\Phi , \Psi ) = (X^{\ast} \Phi , Y \Psi) = (Y^{\ast} \Phi , X \Psi) .
\]
Let $\{ P_{j} \}_{j=1}^{N}$, $N \in \mbf{N}$, be self-adjoint operators on a Hilbert space 
$\ms{H}$.  Set $\ms{D}_{0} = \cap_{j=1}^N \ms{D}(P_{j})$. Assume that $\{ P_{j} \}_{j=1}^{N}$  satisfying the following condition: 
\begin{quote}
\textbf{(S.1)} $\ms{D}_{0}$ is dense in \ms{H}. For all $ \Phi , \Psi \in \ms{D}_{0} $, $ [P_{j} , P_{l}]^0 (\Phi , \Psi ) = 0 $,  $j, l , = 1, \cdots, N$. 
\end{quote}
Let $M$ be a bounded self-adjoint operator $\ms{H}$ satisfying the condition below :
\begin{quote}
 \textbf{(S.2)} For all $\Phi , \Psi \in \ms{D}_{0}$, $ [M, P_{j} ]^0(\Phi , \Psi )=0 $, $j=1 , \cdots , N$. 
\end{quote}
Let $\ms{K}$ be a Hilbert space. Let $\{ \Gamma_{j} \}_{j=1}^{N}$ and $B$ be bounded self-adjoint operators on $\ms{K}$ satisfying the anti-commutation relations below:
\begin{quote}
\textbf{(S.3)}
$
\textbf{(i)}\; \{\Gamma_{j}  , \Gamma_{k}   \}  = 2 \delta_{j,k} , \;  j,k = 1, \cdots , N ,\; \; \textbf{(ii)} \;
\{\Gamma_{j}  , B    \} =1, j= 1, \cdots N , \;  \; \; \; \textbf{(iii)} \; B^2 =I 
$.
\end{quote}
Then the next assertion holds.
\begin{theorem}
Let $ \ms{H}_{\Dir} = \ms{K} \tens \ms{H}$. 
Assume \textbf{(S.1)} - \textbf{(S.3)}. Then 
\[
H_{\Dir}= \sum\limits_{j=1}^{N} \Gamma_{j} \tens P_{j} + B \tens M ,
\]
 is self-adjoint on $\ms{D}(H_{\Dir})=  \ms{K} \tens \ms{D}_{0}$.
\end{theorem}

\begin{remark} In the case where $\{ P_{j}\}_{j=1}^N $ strongly commute, Theorem 2 has been proven in (\cite{Ar93};Theorem 4.3, \cite{Ar07};Lemma 6.7) by strongly anti-commuting methods (\cite{Pe90, Va83}). 
\end{remark}
$ \, $ \\
It is seen that $(I \tens B)^2 =I$ and $(I \tens B)  \ms{D} (H_{\Dir}) \subset \ms{D} (H_{\Dir})$. Then  from Theorem 2,     $H_{\Dir}$ is an abstract Dirac operator on $\ms{H}_{\Dir}$ with the unitary involution $I \tens B$. 

$\; $\\
To prove Theorem 2, we show some lemmas. 
 \begin{lemma} \label{CjCl}
 Let $\{ C_{j} \}_{j=1}^{N}$, $N \in \mbf{N} $, be closed operators on a Hilbert space on $\ms{X}$. Suppose that $ \cap_{j=1}^{N} \ms{D}(C_{j}) $ is dense in $\ms{X}$ and  for $j \ne l$, 
$ ( C_j \Psi , C_l \Psi ) + ( C_l \Psi , C_j \Psi ) = 0 $, $  \Psi \in \cap_{j=1}^{N} \ms{D}(C_{j}) $.
Then $C = \sum\limits_{j=1}^N C_j $ is closed.
\end{lemma}
\textbf{(Proof)}
We see that 
$
(C \Psi , C \Psi ) = \sum\limits_{j=1}^N  \| C_j \Psi \|^2 \geq \frac{1}{N}\left( \sum\limits_{j=1}^N  \| C_j \Psi \| \right)^2 
$.
Then  $ \sum\limits_{j=1}^N \|C_j \Psi \| \leq \sqrt{N} \| C \Psi \|$. Then from a closedness criterion (e.g., \cite{Ar};Theorem B1, \, \cite{GJ69};Proposition 1),  $C$ is closed.  $\blacksquare $

$\; $ \\
From an argument of quadratic forms, there exists a self-adjoint operator $L$ on $\ms{H}$ such that $L \geq 1$, 
$ \ms{D} (\sqrt{L}) =  \ms{D}_{0} $ and for all $\Phi , \Psi \in \ms{D}(\sqrt{L})$,
\begin{equation}
(\sqrt{L}   \Phi , \sqrt{L}   \Psi) \; = \; \; \sum_{j=1}^{N} \left( P_{j} \Phi , P_{j} \Psi  \right) + 
(\Phi , \Psi ).  
\end{equation}

%\begin{lemma} \label{11/23.a} Assume \textbf{(S.1)}-\textbf{(S.3)}. Then $\ms{D} (H_{\Dir}) = \ms{K} \tens \ms{D} (\sqrt{L} )$ and for all $\Psi \in \ms{D} (H_{\Dir}) $ \end{lemma} \textbf{(Proof)}  Since $\ms{D}_0 = \ms{D}(\sqrt{L})$,  $\ms{D} (H_{\Dir}) = \ms{K} \tens \ms{D} (\sqrt{L})$.   It follows that for all $\Psi \in \ms{D} (H_{\Dir})$,  Then we obtain (\ref{norm.1/21.1}).$\blacksquare$

\begin{lemma} \label{rootL}
Assume \textbf{(S.1)}. Then for all $ \Phi , \Psi \in \ms{D}(\sqrt{L})$,
\begin{equation}
\quad \quad  \quad 
[ \frac{}{} \sqrt{L} , P_{j} ]^0 ( \Phi , \Psi  ) = 0 , \quad \quad j= 1, \cdots  N .
\end{equation}
\end{lemma}
\textbf{(Proof)} 
 Since $L$ is  positive and  self-adjoint, it follows that $\sqrt{L} \Xi = \int_{0}^{\infty}\frac{1}{\sqrt{\lambda}} (L+ \lambda)^{-1}L \Xi $, $\Xi \in \ms{D} (L)$. Then for all $ \Phi , \Psi \in \ms{D} (L)$,
  \begin{align}
 [  \sqrt{L} , \, P_{j} ]^0 ( \Phi , \Psi  ) & =  \int_{0}^{\infty}\frac{1}{\sqrt{\lambda}}  
\left[ (L+ \lambda)^{-1} L   , \,  P_{j} \right]^0 (\Phi , \Psi ) d \lambda  \notag \\
& =  \int_{0}^{\infty} \sqrt{\lambda}  
 \left[  P_{j} , \,  L  \right]^0 \left( (L+ \lambda)^{-1} \Phi ,(L+ \lambda)^{-1}  \Psi \right) d \lambda \notag \\
 & = \sum_{l=1}^{N} \int_{0}^{\infty} \sqrt{\lambda}  \left\{
 \left[  P_{j} , \,  P_{l}  \right]^0 \left( (L+ \lambda)^{-1} \Phi ,P_{l} (L+ \lambda)^{-1}  \Psi \right) \right. \notag \\
 & \qquad \quad + \left. \left[  P_{j} , \,  P_{l}  \right]^0 \left( P_{l} (L+ \lambda)^{-1} \Phi , (L+ \lambda)^{-1}  \Psi \right)  \right\} d \lambda  \label{1/23.d} .
 \end{align}
By  \textbf{(S.1)} and (\ref{1/23.d}), we have $[  \sqrt{L} , \, P_{j} ]^0 ( \Phi , \Psi  )=0$ for all $\Phi , \Psi \in \ms{D}(\sqrt{L})$. Note that $\ms{D}(L)$ is a core of $\sqrt{L}$, since $L$ is self-adjoint. In addition, for all $\Psi \in \ms{D}(\sqrt{L})$, $\| P_{j} \Psi \|   \leq \| \sqrt{L} \Psi \|$, $j=1, \cdots , N$. Hence  it follows that  $ [ \sqrt{L} , \, P_{j} ]^0 ( \Phi , \Psi  ) = 0$ for all $\Phi, \Psi \in \ms{D} (\sqrt{L})$. Thus the proof is obtained. $\blacksquare $

$\;$ \\
{\large \textbf{(Proof of Theorem 2)}} \\
Since $B \tens M$ is bounded, it is enough to show that $H= \sum\limits_{j=1}^{N}\Gamma_{j}\tens P_{j} $ is self-adjoint. 
 Let $X = B \tens \sqrt{L}$.  We show that $H$ and $X$ satisfy \textbf{(C.1)}-\textbf{(C.3)} and \textbf{(I)}
  in Theorem 1.    Since $H$ is symmetric  and $X$ self-adjoint, \textbf{(C.1)} is satisfied. 
Since $\sigma (B) = \{ \pm 1\} $ and $ \sqrt{L} \geq 1$, we see that for all $\Psi \in \ms{D}(H)$,
\begin{equation}
| ( \Psi , X \Psi)| =  (\Psi, (I \tens \sqrt{L}) \Psi )  \geq (\Psi , \Psi ).  \label{1/23.b}
\end{equation} 
Then \textbf{(C.2)} is satisfied. Since $\ms{D}_{0}=\ms{D}(\sqrt{L})$, it follows that $\ms{D}(H)= \ms{D}(X)$. Then by \textbf{(S.3)}, we  see that for all $\Psi \in \ms{D}(H)$,
 \begin{equation} 
\| H \Psi \|^2 = \sum_{j=1}^{N} \left( \frac{}{}( I \tens P_{j}) \Psi , (I \tens P_{j}) \Psi  \right) \leq  \| (I \tens \sqrt{L}) \Psi \|^{2} = \| X \Psi \|^2 . \notag
 \end{equation}
Then $\|H \Psi \| \leq \| X \Psi  \|$ for all $\Psi \in \ms{D}(H)$, and hence \textbf{(C.3)} is satisfied.
By   Lemma \ref{rootL}, it is seen that for all $\Psi \in \ms{D} (H)$, 
\begin{align}
( H \Psi , \, X \Psi ) + ( X \Psi , H \Psi ) & = \sum_{j=1}^N
\left(   ( (\Gamma_{j}\tens P_{j} ) \Psi , ( B \tens \sqrt{L}) \Psi ) 
+ (  (B \tens \sqrt{L}) \Psi  ,( \Gamma_{j}\tens P_{j}) \Psi    )   \right) \notag \\
& = \sum_{j=1}^N
\left(   ( I \tens \sqrt{L} ) \Psi , (\{ \Gamma_{j} ,  B \} \tens P_{j})  \Psi )    \right)   \label{1/23.c}
\end{align}
Then by \textbf{(S.3)} and (\ref{1/23.c}), we have $( H \Psi , \, X \Psi ) + ( X \Psi , H \Psi ) =0$. Then from (\ref{1/23.b}), it follows that $( H \Psi , \, X \Psi ) + ( X \Psi , H \Psi ) \leq | ( \Psi ,X \Psi )| $ 
for all $ \Psi \in \ms{D} (H)$. Then \textbf{(I)} is satisfied, and hence  $\overline{H}$  is  self-adjoint  from Theorem 1. In addition, by \textbf{(S.1)} and \textbf{(S.3)}, we see that  for $j \ne l$,
\[
 (  (\Gamma_{j}\tens P_{j}) \Psi , ( \Gamma_{l}\tens P_{l})  \Psi  ) + ( ( \Gamma_{l}\tens P_{l}) \Psi ,
( \Gamma_{j}\tens P_{j} ) \Psi  )
 = (( I \tens P_{l}) \Psi ,( \{ \Gamma_{j}, \Gamma_{l} \} \tens P_{j}) \Psi  )
 = 0 .
 \]
Then from Lemma \ref{CjCl}, $\overline{H} = H$, and hence the proof is obtained. $\blacksquare$

$\;$ \\
{\large \textbf{Acknowledgments }} 
It is a pleasure to thank  Professor Akito Suzuki and Professor  Fumio Hiroshima for their  comments.
This work is supported by JSPS grant  24$\cdot$1671.

\end{document}